\documentclass[notitlepage]{amsart}

\usepackage{amsfonts}
\usepackage[latin1]{inputenc}
\usepackage[all]{xy} 
\usepackage{latexsym,amssymb,amscd}
\usepackage{amsmath}
\newtheorem{pro}{Proposition}[section]
\newtheorem{teo}[pro]{Theorem}
\newtheorem{defi}[pro]{Definition}
\newtheorem{lem}[pro]{Lemma}
\newtheorem{cor}[pro]{Corollary}
\newtheorem{rk}[pro]{Remark}

\newcommand{\la}{\left\langle}
\newcommand{\ra}{\right\rangle}
\newcommand{\Ext}{\mathrm{Ext}}
\newcommand{\Hom}{\mathrm{Hom}}

\newcommand{\U}{\mathcal{U}}

\newcommand{\T}{\mathcal{T}}

\newcommand{\X}{\mathcal{X}}
\newcommand{\Y}{\mathcal{Y}}
\newcommand{\Z}{\mathcal{Z}}

\newcommand{\N}{\mathbb{N}}
\newcommand{\Enteros}{\mathbb{Z}}
\newcommand{\pd}{\mathrm{pd}}
\newcommand{\id}{\mathrm{id}}
\newcommand{\resdim}{\mathrm{resdim}}
\newcommand{\coresdim}{\mathrm{coresdim}}
\newcommand{\add}{\mathrm{add}}
\newcommand{\mini}{\mathrm{min}}
\newcommand{\maxi}{\mathrm{max}}

\newenvironment{dem}{\noindent\bf Proof. \rm }{$\ \Box$}
\usepackage{latexsym,amssymb,amscd}
\usepackage{amsmath}

\newcommand{\Qed}{~\hfill$\square$\bigskip}

\begin{document}
\title[Auslander-Buchweitz approximation theory]{Auslander-Buchweitz approximation theory for triangulated categories}
\author{O. Mendoza, E. C. S\'aenz, V. Santiago, M. J. Souto Salorio.}
\thanks{
The authors thank the financial support received from
Project PAPIIT-UNAM IN101607 and MICINN-FEDER   TIN2010-18552-C03-02.}
\date{}
\begin{abstract}
We introduce and develop an analogous of the Auslander-Buchweitz
approximation theory (see \cite{AB}) in the context of triangulated
categories, by using a version of relative homology in this setting.
We also prove several results concerning relative homological
algebra in a triangulated category $\T,$ which are based on the
behaviour of certain subcategories under finiteness of resolutions
and vanishing of Hom-spaces. For example: we establish the existence
of preenvelopes (and precovers) in certain triangulated
subcategories of $\T.$ The results resemble various constructions
and results of Auslander and Buchweitz, and are concentrated in
exploring the structure of a triangulated category $\T$ equipped
with a pair $(\X,\omega),$ where $\X$ is closed under extensions and
$\omega$ is a weak-cogenerator in $\X,$ usually under additional
conditions. This reduces, among other things, to the existence of
distinguished triangles enjoying special properties, and the
behaviour of (suitably defined) (co)resolutions, projective or
injective dimension of objects of $\T$ and the formation of
orthogonal subcategories. Finally, some relationships with the Rouquier's dimension
in triangulated categories is discussed.
\end{abstract}
\maketitle
\section*{Introduction.}

The approximation theory has its origin with the concept of injective envelopes and it has had a wide development in the context of module categories since the fifties.

In independent papers,  Auslander, Reiten and Smalo (for the category $\mathrm{mod}\,(\Lambda)$ of finitely generated modules over an artin algebra $\Lambda$), and Enochs (for the category $\mathrm{Mod}\,(R)$ of modules over an arbitrary ring $R$) introduced a general approximation theory involving precovers and preenvelopes
(see \cite{AR}, \cite{AS} and \cite{E}).

Auslander and Buchweitz (see \cite{AB})
studied the ideas of injective envelopes and projective covers
in terms of maximal Cohen-Macaulay
approximations for certain modules. In their work, they also studied
the relationship between the relative injective dimension and
the coresolution dimension of a module. They
developed their theory in the context of
abelian categories providing important applications in several settings.

 Based on \cite{AB}, Hashimoto defined the so called ``Auslander-Buchweitz context" for abelian categories, giving a new framework to homological approximation theory  (see \cite{H}).

Recently, triangulated categories entered into the subject in a relevant way and several authors have studied the concept of approximation  in both contexts, abelian and triangulated categories (see, for example, \cite{AM}
\cite{B}, \cite{BR} and \cite{MS}).

In this paper, an analogous theory of approximations in the sense of Auslander and Buchweitz
(see \cite{AB}), is developed for triangulated categories. Throughout this paper, $\T$ denotes 
an arbitrary triangulated category and $\X$ a class of objects in $\T$.
The main result (Theorem \ref{specialtrian}) deals with a pair $(\X,\omega)$ of classes of 
objects in $\T ,$ where $\X$ is closed under extensions, and $\omega$ satisfies a weak 
cogenerating condition with respect to the objects of $\X.$ Like Auslander and Buchweitz, we
consider the class $\X^\wedge$ of objects of $\T$ admitting  a finite resolution by objects of 
$\X.$ We prove that any object of $\X^\wedge$ admits two distinguished triangles:
one giving rise to a right $\X$-approximation, and the other to a left 
$\omega^\wedge$-approximation. In the present paper, it is also introduced and discussed a notion 
of $\X$-resolution dimension, which is compared with other relative homological dimensions.

The paper is organized as follows:
In Section 1, we give some basic notions and properties of triangulated categories, that will 
be used in the rest of the work.

In Section 2, we study the notion of $\mathcal{X}$-resolution
dimension which allows us to characterize the triangulated subcategory $\Delta_\T(\X)$ of $\T,$ generated by a cosuspended subcategory $\X$ of $\T$ (see Theorem \ref{resBuan}).

In Section 3, the properties of the  $\mathcal{X}$-projective (respectively,
$\mathcal{X}$-injective) dimension and its relation to the
$\mathcal{X}$-resolution (respectively, coresolution) dimension are established.
The main result of this section is Theorem \ref{pdyresdim} that relates different
kinds of relative homological dimensions by using suitable subcategories
of $\mathcal{T}$.

In Section 4, we focus our attention to the notions of
$\mathcal{X}$-injectives and weak-cogenerators in $\mathcal{X}.$ We
relate these ideas to the concepts of injective and coresolution
dimension. This leads us to characterize several triangulated
subcategories; and moreover, in Theorem \ref{specialtrian} we
establish the existence of $\mathcal{X}$-precovers and
$\omega^\wedge$-preenvelopes. Finally, in Theorem
\ref{tilde-id}, for  a given pair $(\X,\omega)$ satisfying certain
conditions, we give several characterizations of the triangulated
subcategory $\Delta_{\X^\wedge}(\omega)$ of $\X^\wedge$ generated by
$\omega.$

We remark that the results we get  will be applied in a forthcoming paper
\cite{MSSS2}, where a connection between Auslander-Buchweitz approximation theory in 
triangulated categories and co-$t$-structures (see \cite{Bo} and \cite{P}) is established.

Finally, some relationship with other notions, as torsion theories (see Corollary \ref{torsion}) in the sense of Iyama-Yoshino 
\cite{IY} and Rouquier's dimension \cite{Ro}, are discussed (see Section 5).

\section{Preliminaries}
\

Throughout this paper, $\T$ will be a triangulated category and
$[1]:\T\rightarrow\T$ its suspension functor. The term subcategory, in this paper, means a subcategory which is full,
additive, and closed under isomorphisms.
\

An important tool, which is a consequence of the octahedral axiom in $\T,$ is the so-called co-base change (see \cite{K}). That is, for any diagram in $\T$
$$\begin{CD}
X @>>> Y\\
@VVV @.\\
Z @.\,
\end{CD}$$
\vspace{.2cm}
 there exists a commutative and exact diagram in $\T$

$$\begin{CD}
@. W[-1] @= W[-1] @.\\
@. @VVV @VVV @.\\
U[-1] @>>> X @>>> Y @>>> U\\
@| @VVV @VVV @|\\
U[-1] @>>> Z @>>> E @>>> U\\
@. @VVV @VVV @.\\
@. W @= W @.\,  \\
{}
  \end{CD}
  $$
where exact means that the rows and columns, in the preceding diagram, are distinguished triangles in $\T.$ The base change, which is the dual notion of co-base change, also holds (see \cite{K}).

Let $\X$ and $\Y$ be classes of objects in $\T.$ We put
${}^\perp\X:=\{Z\in\T\,:\,\Hom_\T(Z,-)|_{\X}=0\}$ and $\X^\perp:=\{Z\in\T\,:\,\Hom_\T(-,Z)|_{\X}=0\}.$
\noindent We denote by $\X*\Y$ the class of objects $Z\in\T$ for which exists a distinguished 
triangle $X\rightarrow Z\rightarrow Y\rightarrow X[1]$ in $\T$ with $X\in\X$ and $Y\in\Y.$ In case 
$\Y=\{Y\},$ we write $\X*Y$ instead of $\X*\Y.$\\
It is also well known that the operation $*$ is associative (see \cite[1.3.10]{BBD}). Furthermore, it is said that  $\X$ is {\bf closed under extensions} if $ \X*\X\subseteq \X.$
\

Recall that  a class $\X$ of objects in $\T$ is said to be {\bf
suspended} (respectively,  {\bf cosuspended}) if $\X[1]\subseteq\X$
(respectively,  $\X[-1]\subseteq\X$) and $\X$ is closed under
extensions. By the following lemma, it is easy to see, that a
suspended (respectively, cosuspended) class $\X$ of objects in $\T,$
can be considered as a subcategory of $\T.$

\begin{lem}\label{xx=x}
Let $\X$ be a class of objects in  $\T.$
\begin{enumerate}
\item[(a)]  If $0\in \X$  then $\Y\subseteq\X*\Y$ and $\Y\subseteq\Y*\X$ for any class $\Y$ of objects in  $\T.$
\item[(b)] If $\X$ is either suspended or cosuspended, then $0\in \X$  and  $\X=\X*\X.$
\end{enumerate}
\end{lem}
\begin{dem} (a) If $0\in \X$  then we get $\Y\subseteq \X*\Y$  by using the distinguished triangle $0\rightarrow Y\stackrel{1_Y}{\to}Y \to 0$ for any $Y\in \Y.$ The other inclusion follows similarly.
\

(b) Let $\X$ be cosuspended (the other case, is analogous). Then, it follows that $0\in \X$ since we have the distinguished triangle $X[-1]\rightarrow 0\rightarrow X  \rightarrow X $ for any  $X\in \X.$ Hence (b) follows from (a).
\end{dem}
\vspace{.2cm}

Given a class $\X$ of objects in $\T,$ it is said that $\X$ is {\bf{closed under cones}} if for any distinguished
triangle  $A\rightarrow B\rightarrow C\rightarrow A[1]$ in $\T$ with  $A,B\in \X$ we have that $C\in \X.$  Similarly,
$\X$ is {\bf{closed under cocones}} if for any distinguished triangle  $A\rightarrow B\rightarrow C\rightarrow A[1]$ in
$\T$ with $B,C\in \X$ we have that $A\in \X.$
\

We denote by $\U_{\X}$ (respectively,  ${}_{\X}\U$) the smallest
suspended (respectively,  cosuspended) subcategory of $\T$
containing the class $\X.$ Note that if $\X$ is suspended
(respectively, cosuspended) subcategory of $\T,$ then $\X=\U_{\X}$
(respectively,  $\X={}_{\X}\U$). We also recall that a subcategory
$\U$ of $\T,$ which is suspended and cosuspended, is called
{\bf{triangulated subcategory}} of $\T.$ A {\bf{thick}} subcategory
of $\T$ is a triangulated subcategory of $\T$ which is closed under
direct summands in $\T.$ We also denote by $\Delta_\T(\X)$
(respectively, $\overline{\Delta}_\T(\X))$ to the smallest triangulated
(respectively, smallest thick) subcategory of $\T$ containing the
class $\X.$ Observe that
$\Delta_\T(\X)\subseteq\overline{\Delta}_\T(\X).$ For the following definition, see \cite{AR}, \cite{B}, \cite{BR} and \cite{E}.

\begin{defi} Let $\X$ and $\Y$ be classes of objects in the triangulated category $\T.$ A morphism $f:X\to C$ in $\T$ is said to be an {\bf $\X$-precover} of $C$ if $X\in\X$ and $\Hom_\T(X',f):\Hom_\T(X',X)\to\Hom_\T(X',C)$ is surjective, $\forall X'\in\X.$ If any $C\in\Y$ admits an $\X$-precover, then $\X$ is called a {\bf precovering class} in $\Y.$ By dualizing the definition above, we get the notion of an {\bf $\X$-preenveloping} of $C$ and a {\bf preenveloping class} in $\Y.$
\end{defi}
\

Finally, in order to deal with the (co)resolution, relative projective and relative injective dimensions, we consider the extended
natural numbers $\overline{\mathbb{N}}:=\mathbb{N}\cup\{\infty\}.$ Here, we set the following rules: (a) $x+\infty=\infty$ for any
$x\in\overline{\mathbb{N}},$ (b) $x<\infty$ for any $x\in\mathbb{N}$ and (c) $\mini(\emptyset):=\infty.$

\section{resolution and coresolution dimensions}

Now, we define certain classes of objects in $\T$ which will lead us to the notions of resolution and coresolution dimensions.

\begin{defi}\label{epsilon} Let $\X$ be a class of objects in $\T.$ For any natural number $n,$ we introduce inductively the  class $\varepsilon^\wedge_n(\X)$ as follows: $\varepsilon^\wedge_0(\X):=\X$ and assuming defined $\varepsilon^\wedge_{n-1}(\X),$ the class $\varepsilon^\wedge_{n}(\X)$ is given by all the objects $Z\in\T$ for which there exists a distinguished triangle in $\T$
$$\begin{CD}
Z[-1] @>>> W @>>> X @>>> Z \,
\end{CD}$$ with $W\in\varepsilon^\wedge_{n-1}(\X)$ and $X\in\X.$

Dually, we set
 $\varepsilon^\vee_0(\X):=\X$ and supposing defined $\varepsilon^\vee_{n-1}(\X),$ the class $\varepsilon^\vee_{n}(\X)$ is formed for all the objects $Z\in\T$ for which there exists a distinguished triangle in $\T$
$$\begin{CD}
Z @>>> X @>>> K @>>>  Z[1]\,
\end{CD}$$ with $K\in\varepsilon^\vee_{n-1}(\X)$ and $X\in\X.$
\end{defi}

We have the following properties for $\varepsilon^\wedge_{n}(\X)$
(and the similar ones for $\varepsilon^\vee_{n}(\X)$).

\begin{pro}\label{rkepsilon} Let $\T$ be a triangulated category, $\X$ be a class of objects in $\T,$ and $n$ a natural number. Then, the following statements
hold.
\begin{enumerate}
\item[(a)] For any $Z\in \T$ and $n>0,$ we have that $Z\in\varepsilon^\wedge_{n}(\X)$ if and only if there is a family
$\{K_j[-1]\to K_{j+1}\to X_j\to K_j\}_{j=0}^{n-1}$ of
distinguished triangles in $\T$ with $K_0=Z,$ $X_j\in \X$ and
$K_n\in\X.$
\item[(b)] $\varepsilon^\wedge_{n}(\X)=*_{i=0}^n\,\X[i]:=\X*\X[1]*\cdots*\X[n].$
\item[(c)] If $0\in \X$ then $\X[n]\subseteq
\varepsilon^\wedge_{n}(\X)\subseteq\varepsilon^\wedge_{n+1}(\X)$ and
$\varepsilon^\wedge_n(\X)[1]\subseteq \varepsilon^\wedge_{n+1}(\X)$ $\quad\forall\;n\in\N.$
\end{enumerate}
\end{pro}
\begin{dem} (a) If $n=1$ then the equivalence follows from the definition of $\varepsilon^\wedge_1(\X).$ Let $n\geq 2$
and suppose (by induction) that the equivalence is true for
$\varepsilon^\wedge_{n-1}(\X).$ By definition, $Z\in
\varepsilon^\wedge_{n}(\X)$ if and only if there is a distinguished
triangle in $\T$
$$\begin{CD}
Z[-1] @>>> K_1 @>>> X_0 @>>> Z \,
\end{CD}$$ with $K_1\in\varepsilon^\wedge_{n-1}(\X)$ and $X_0\in\X.$ On the other hand, by induction, we have that
$K_1\in\varepsilon^\wedge_{n-1}(\X)$ if and only if there is a
family $\{K_j[-1]\to K_{j+1}\to X_j\to K_j\}_{j=1}^{n-1}$
of distinguished triangles in $\T$ with $X_j\in \X$ and $K_n\in\X;$
proving (a).
\

(b) By definition, we have that
$\varepsilon^\wedge_{n}(\X)=\X*\varepsilon^\wedge_{n-1}(\X)[1].$ So,
by induction, it follows that
$\varepsilon^\wedge_{n}(\X)=\X*(*_{i=0}^{n-1}\X[i])[1]=*_{i=0}^n\,\X[i].$
\

(c) Assume that $0\in \X.$ By (b), we know that 
$\varepsilon^\wedge_{n}(\X)=\varepsilon^\wedge_{n-1}(\X)*\X[n];$ and since 
$0\in\varepsilon^\wedge_{n-1}(\X),$ it follows from \ref{xx=x} (a) that 
$\X[n]\subseteq\varepsilon^\wedge_{n}(\X).$ Similarly, from the equalities 
$\varepsilon^\wedge_{n+1}(\X)=\varepsilon^\wedge_n(\X)*\X[n+1]$ and 
$\varepsilon^\wedge_{n+1}(\X)=\X*\varepsilon^\wedge_n(\X)[1],$ and the facts that
$0\in \X[n+1]$ and $0\in\X,$ we get the other inclusions from \ref{xx=x} (a).
\end{dem}
\vspace{.2cm}

Following \cite{AB} and \cite{B}, we introduce the notion of $\X$-resolution (respectively, coresolution) dimension of any class $\Y$ of objects of $\T.$

\begin{defi} Let $\X$ be a class of objects in $\T.$
\begin{enumerate}
 \item[(a)] $\X^\wedge:=\cup_{n\geq 0}\;\varepsilon^\wedge_n(\X)$  and $\X^\vee:=\cup_{n\geq 0}\;\varepsilon^\vee_n(\X).$
 \item[(b)] For any $M\in\T,$ the $\X$-{\bf{resolution dimension}} of $M$ is
$$\resdim_{\X}(M):=\mini\,\{n\in\mathbb{N} \;:\; M\in\varepsilon^\wedge_{n}(\X) \}.$$ Dually, the $\X$-{\bf{coresolution dimension}} of $M$ is
$$\coresdim_{\X}(M):=\mini\,\{n\in\mathbb{N}:\; M\in\varepsilon^\vee_{n}(\X) \}.$$
 \item[(c)] For any subclass $\Y$ of $\T$, we set $\resdim_{\X}(\Y):=\sup\,\{\resdim_{\X}(M)\,:\,M\in\Y\}.$ Similarly, we also have $\coresdim_{\X}(\Y).$
\end{enumerate}
\end{defi}

\begin{pro}\label{stardim} Let $\X$ and $\Y$ be   classes of objects in $\T,$ $0\in\Y$ and $n\in\N.$ Then, the following statements hold.
 \begin{itemize}
  \item[(a)] $\resdim_\Y(\X)\leq n$ if and only if $\X\subseteq\varepsilon^\wedge_{n}(\Y)= *_{i=0}^n\,\Y[i].$
  \item[(b)] If $\X$ is closed under extensions, then  $\X* \X^{\wedge} \subseteq  \X^{\wedge}.$
  \item[(c)] If $\X$ is cosuspended, then $\varepsilon^\wedge_{n}(\X)=\X[n].$ 
 \end{itemize}
\end{pro}
\begin{dem} (a) It follows by definition and \ref{rkepsilon} (b),(c).
\

(b) It follows from \ref{rkepsilon} (b) since $\X*\X\subseteq\X.$
\

(c) Let $\X$ be cosuspended. since $0\in\X$ (see \ref{xx=x} (b)),
we get from \ref{rkepsilon} (b), (c) that $\X[n]\subseteq\varepsilon^\wedge_{n}(\X)= *_{i=0}^n\,\X[i].$ On the other hand, using that $\X*\X\subseteq\X$
and $\X[-1]\subseteq \X,$ we conclude that
$*_{i=0}^n\,\X[i]=(*_{i=0}^n\,\X[i-n])[n]\subseteq(*_{i=0}^n\,\X)[n]\subseteq\X[n].$
\end{dem}
\vspace{.2cm}

The following result will be useful in this paper. The item (a) already appeared in \cite{B}. We also recall that $\Delta_\T(\X)$
(respectively, $\overline{\Delta}_\T(\X))$ stands for
the smallest triangulated (respectively, smallest thick) subcategory of $\T$ containing the class of objects $\X.$

\begin{teo}\label{resBuan} For any cosuspended subcategory $\X$ of $\T$ and any object $C\in\T,$  the following statements hold.
 \begin{itemize}
  \item[(a)] $\resdim_\X(C)\leq n$ if and only if $C\in\X[n].$
  \item[(b)] $\X^\wedge=\cup_{n\geq 0}\,\X[n]=\Delta_\T(\X).$
  \item[(c)] If $\X$ is closed under direct summands in $\T,$ then $\X^\wedge=\overline{\Delta}_\T(\X).$
 \end{itemize}
\end{teo}
\begin{dem} (a) If follows from \ref{stardim} (a), (c) since $0\in\X$ (see \ref{xx=x} (b)).
\

 (b) From \ref{stardim} (c), we get $\X^\wedge=\cup_{n\geq 0}\,\X[n];$ and hence $\X^\wedge$ is closed under positive and negative shifts. We prove now
that $\X^\wedge$ is closed under extensions. Indeed, let $X[n]\to Y\to X'[m]\to X[n][1]$ be a distinguished triangle in $\T$
with $X,X'\in\X.$ We assume that $n\leq m$ and then $X[n]=X[n-m][m]\in\X[m]$ since $n-m\leq 0$ and  $\X[-1]\subseteq\X.$ Using now
that $\X$ is closed under extensions, it follows that $Y\in\X[m]\subseteq \X^\wedge;$ proving that $\X^\wedge$ is closed under extensions.
Hence $\X^\wedge$ is a triangulated subcategory of $\T$ and moreover it is the smallest one containing $\X$ since
$\X^\wedge=\cup_{n\geq 0}\,\X[n].$
\

(c) It follows from (b).
\end{dem}

\begin{rk}\label{varios}
$(1)$ Observe that a suspended class $\,\U$ of $\T$ is closed under cones. Indeed, if $A\rightarrow B\rightarrow C\rightarrow A[1]$ is a distinguished triangle in $\T$ with  $A,B\in \U$  then $A[1],B\in \U;$ and so we get  $C\in \U.$ Similarly, if $\,\U$ is cosuspended then it is closed under cocones.
\

$(2)$ Let $(\Y,\omega)$ be a pair of classes of objects in $\T$ with $\omega\subseteq\Y.$ If $\Y$ is closed under cones (respectively, cocones) then $\omega^\wedge\subseteq\Y$ (respectively, $\omega^\vee\subseteq\Y$). Indeed, assume that $\Y$ is closed under cones and let $M\in\omega^\wedge.$ Thus $M\in\varepsilon^\wedge_n(\omega)$ for some $n\in\Bbb{N}.$ If $n=0$ then $M\in\omega\subseteq \Y.$ Let $n>0,$ and hence there is a distinguished triangle $M[-1] \to K\to Y\to M$ in $\T$ with $K\in\varepsilon^\wedge_{n-1}(\omega)$ and $Y\in\Y.$ By induction $K\in\Y$ and hence $M\in\Y$ since $\Y$ is closed under cones; proving that $\omega^\wedge\subseteq\Y.$
\

$(3)$ Note that $\X^\wedge\subseteq  \U_{\X}$ (respectively,  $\X^\vee \subseteq {}_{\X}\U$) since
$\U_{\X}$ (respectively,  ${}_{\X}\U$)  is closed under cones (respectively,  cocones)  and contains $\X.$
\end{rk}

 Using the fact that the functor $\Hom$ is a cohomological one, we get the following description of the orthogonal categories. In particular, observe that ${}_{\X}\U^\perp$ (respectively,  ${}^\perp\U_{\X}$) is a suspended (respectively,  cosuspended) subcategory of $\T.$

\begin{lem}\label{susp-perp} For any class $\X$ of objects in $\T,$ we have that
 \begin{itemize}
 \item[(a)] ${}^\perp\U_{\X}=\{Z\in\T\;:\;\Hom_\T(Z,X[i])=0,\quad \forall i\geq 0, \forall X\in\X\},$
 \item[(b)] ${}_{\X}\U^\perp=\{Z\in\T\;:\;\Hom_\T(X[i],Z)=0,\quad \forall i\leq 0, \forall X\in\X\}.$
 \end{itemize}
\end{lem}
\begin{dem} It is straightforward.
\end{dem}

\begin{lem}\label{restri} Let $\Y$ and $\X$ be   classes  of objects in $\T,$  $n\geq 1$ and $Z\in \T.$
The following statements hold.
 \begin{itemize}
  \item[(a)] The object $Z$ belongs to $\Y*\Y[1]*\cdots*\Y[n-1]* \X[n]$ if and only if there exists a family
   $\{ K_i\rightarrow Y_i\rightarrow K_{i+1}\rightarrow K_i[1]\;:\; Y_i\in\Y\}_{i=0}^{n-1}$ of distinguished triangles in $\T$ with $K_0\in \X$ and $Z=K_n.$
\vspace{.2cm}
  \item[(b)] The object $Z$ belongs to $\X[-n]  *\Y[-n+1]*\cdots*\Y[-1]*\Y$ if and only if there exists a family
   $\{K_{i+1}\rightarrow Y_i\rightarrow K_{i}\rightarrow K_{i+1}[1]\;:\; Y_i\in\Y\}_{i=0}^{n-1}$ of distinguished triangles in $\T$ with $K_0\in \X$ and  $Z=K_n.$
 \end{itemize}
\end{lem}
\begin{dem} (a) We proceed by induction on $n.$  If $n=1$ then (a) is trivial. Suppose that $n\geq 2$ and consider the class $$\Z_{n-1}:=\Y*\Y[1]*\cdots*\Y[n-2] * \X[n-1].$$ It is clear that $\Y*\Y[1]*\cdots*\Y[n-1]* \X[n]=\Y*\Z_{n-1}[1];$ and then, we have that $Z\in\Y*\Y[1]*\cdots*\Y[n-1]* \X[n]$ if and only if there is a distinguished triangle
$$\begin{CD} K  @>>> Y  @>>>   Z @>>> K[1]\end{CD}$$ in $\T$ with $Y\in\Y$ and $K\in\Z_{n-1}.$ On the other hand, by induction,
 we have that $K\in\Z_{n-1}$ if and only if there is a family
   $\{ K_i\rightarrow Y_i\rightarrow K_{i+1}\rightarrow K_i[1]\;:\; Y_i\in\Y\}_{i=0}^{n-2}$ of distinguished triangles in $\T$ with $K_0\in \X$ and $K=K_{n-1}.$ So the result follows by adding the triangle above to the preceding family of triangles.
\

(b) It is similar to (a).
\end{dem}


\section{Relative homological dimensions}

In this section, we introduce the $\X$-projective (respectively, injective) dimension of objects in $\T.$ Moreover, we establish a result that relates this relative projective dimension with the resolution dimension as can be seen in Theorem \ref{pdyresdim}.

\begin{defi} Let $\X$ be a class of objects in $\T$ and $M$ an object in $\T.$
\begin{enumerate}
 \item[(a)] The {\bf{$\X$-projective dimension}} of
$M$ is
$$\pd_{\X}(M):=\mini\,\{n\in\mathbb{N}\; :\; \Hom_\T(M[-i],-)\mid_{\X}=0, \quad\forall
i>n\}.$$
 \item[(b)] The {\bf{$\X$-injective dimension}} of $M$ is
$$\id_{\X}(M):=\mini\,\{n\in\mathbb{N}\; :\; \Hom_\T(-,M[i])\mid_{\X}=0, \quad\forall i>n\}.$$
 \item[(c)] For any class $\Y$ of objects in $\T,$ we set
$$\pd_{\X}(\Y):=\sup\,\{\pd_{\X}(C)\,:\,C\in\Y\}\quad\text{and}\quad\id_{\X}(\Y):=\sup\,\{\id_{\X}(C)\,:\,C\in\Y\}.$$
\end{enumerate}
\end{defi}

\begin{lem}\label{compdim} Let $\X$ be a class of objects in $\T.$ Then, the following statements hold.
 \begin{itemize}
 \item[(a)] For any $M\in\T$ and $n\in\Bbb{N},$ we have that
  \begin{itemize}
   \item[(a1)] $\pd_{\X}(M)\leq n$ if and only if $M\in{}^\perp\U_{\X}[n+1];$
   \item[(a2)] $\id_{\X}(M)\leq n$ if and only if $M\in{}_{\X}\U^\perp[-n-1].$
  \end{itemize}
 \item[(b)] $\pd_{\Y}(\X)=\id_{\X}(\Y)$ for any class $\Y$ of objects in $\T.$
 \end{itemize}
\end{lem}
\begin{dem} (a) follows from \ref{susp-perp}, and (b) is straightforward.
\end{dem}

\begin{pro}\label{pdXresdim} Let $\X$ be a class of objects in $\T$ and $M\in\T.$ Then
$$\pd_{\X}(M)=\resdim_{{}^\perp\U_{\X}[1]}(M)\;\text{ and }\;\id_{\X}(M)=\coresdim_{{}_{\X}\U^\perp[-1]}(M).$$
\end{pro}
\begin{dem} Since ${}^\perp\U_{\X}$ is cosuspended (see \ref{susp-perp} (a)), the first equality follows from \ref{compdim} (a1) and \ref{resBuan} (a). The second equality can be proven similarly.
\end{dem}
\vspace{.2cm}

Now, we prove the following relationship between the relative projective dimension and the resolution dimension. 

\begin{teo}\label{pdyresdim} Let $\X$ and $\Y$ be classes of objects in $\T.$ Then, the following statements hold.
\begin{enumerate}
 \item[(a)] $\pd_{\X}(L)\leq \pd_{\X}(\Y) + \resdim_{\Y}(L),$ $\quad \forall L\in\T.$
\vspace{.2cm}
\item[(b)] If $\Y\subseteq \U_{\X}\cap{}^{\perp}\U_{\X}[1]$ and $\Y$ is closed under
direct summands in $\T,$ then $$\pd_{\X}(L)=\resdim_{\Y}(L), \quad  \forall L\in\Y^\wedge.$$
\end{enumerate}
\end{teo}
\begin{dem} (a) Let $d:=\resdim_{\Y}(L)$ and $\alpha:=\pd_{\X}(\Y).$ We may assume that $d$ and $\alpha$ are finite.
We prove (a) by induction on $d.$ If $d=0,$ it follows that $L\in\Y$ and then (a) holds in this case.\\
Assume that $d\geq 1.$ So we have a distinguished triangle $K\rightarrow Y\rightarrow L\rightarrow K[1]$ in $\T$ with $Y\in\Y$ and $K\in\varepsilon^\wedge_{d-1}(\Y).$ Applying the cohomological functor $\Hom_\T(-,M[j]),$ with $M\in\X,$ to the above triangle, we get and exact sequence of abelian groups
$$\Hom_\T(K[1],M[j])\rightarrow \Hom_\T(L,M[j])\rightarrow \Hom_\T(Y,M[j]).$$ By induction, we know that $\pd_\X(K)\leq \alpha+d-1.$ Therefore $\Hom_\T(L,M[j])=0$ for $j>\alpha+d$ and so $\pd_\X(L)\leq \alpha+d.$
\

(b) Let $\Y\subseteq \U_{\X}\cap ^{\perp}\U_{\X}[1]$ and $\Y$ be closed under direct summands in $\T.$ Consider $L\in\Y^\wedge$ and let $d:=\resdim_{\Y}(L).$ By \ref{compdim} we have that $\pd_{\X}(\Y)=0$ and then $\pd_{\X}(L)\leq d$ (see (a)). We prove, by induction on $d,$ that the equality given in (b) holds. For $d=0$ it is clear.\\
Suppose that $d=1.$ Then, there is a distinguished triangle $$\eta:\quad Y_1\rightarrow Y_0\rightarrow L\stackrel{f}{\rightarrow} Y_1[1]\text{ in $\T$ with $Y_i\in\Y.$}$$ If $\pd_{\X}(L)=0$ then $L\in{}^{\perp}\U_{\X}[1]$ (see \ref{compdim}). Hence $f=0$ since $\Y\subseteq \U_{\X};$ and therefore $\eta$ splits giving us that $L\in\Y,$ which is a contradiction since $d=1.$ So  $\pd_{\X}(L)>0$ proving (b) for $d=1.$\\
Assume now that $d\geq 2.$ Thus we have a distinguished triangle $K\rightarrow Y\rightarrow L\rightarrow K[1]$ in $\T$ with $Y\in\Y,$ $K\in\varepsilon^\wedge_{d-1}(\Y)$ and $\pd_\X(K)=d-1$ (by inductive hypothesis). Since $\pd_{\X}(L)\leq d,$ it is enough to see $\pd_{\X}(L)> d-1.$ So, in case $\pd_{\X}(L)\leq d-1,$ we apply the cohomological functor $\Hom_\T(-,X[d]),$ with $X\in\X,$ to the triangle $L\rightarrow K[1]\rightarrow Y[1]\rightarrow L[1].$ Then we get the following exact sequence of abelian groups
$$\Hom_\T(Y[1],X[d])\rightarrow \Hom_\T(K[1],X[d])\rightarrow \Hom_\T(L,X[d]).$$
\vspace{.2cm}
\noindent Therefore $\Hom_\T(K[1],X[d])=0$ contradicting that $\pd_\X(K)=d-1.$ This means that $\pd_{\X}(L)> d-1;$ proving (b).
\end{dem}

\begin{rk} Note that if $Y\neq 0$ and $Y\in \U_{\X}\cap{}^{\perp}\U_{\X}[1],$ then $Y[j]\notin \U_{\X}\cap{}^{\perp}\U_{\X}[1],$ $\forall j>0.$
\end{rk}

The following technical result will be used in Section 4.

\begin{lem}\label{pdim} Let $\X$, $\Y$ and $\Z$ be classes of objects in $\T.$ Then, the following statements hold.
 \begin{enumerate}
  \item[(a)] $\pd_{\Y}(\X^\vee)=\pd_{\Y}(\X).$
\vspace{.2cm}
  \item[(b)] If $\X\subseteq \Z\subseteq \X^\vee$ then
$\pd_{\Y}(\Z)=\pd_{\Y}(\X).$
\end{enumerate}
\end{lem}
\begin{dem} To prove (a), it is enough to see that $\pd_{\Y}\,(\X^\vee) \leq\pd_{\Y}\,(\X).$ Let $M\in\X^\vee$. We prove by induction on $d:=\coresdim_{\X}\,(M)$ that $\pd_{\Y}\,(M) \leq\pd_{\Y}\,(\X)$. We may assume that $\alpha:=\pd_{\Y}\,(\X)<\infty.$ If $d=0$ then we have that $M\in\X$ and there is nothing to prove.
\

Let $d\geq 1.$ Then we have a distinguished triangle $M\rightarrow X\rightarrow K\rightarrow M[1]$ in $\T$ with $X\in\X,$ $K\in\varepsilon^\vee_{d-1}(\X)$ and $\pd_{\Y}\,(K)\leq \alpha$ (by inductive hypothesis). Applying the cohomological functor $\Hom_\T(-,Y[i])$, with $Y\in\Y$, we get the exact sequence of abelian groups
$$\Hom_\T(X,Y[i])\rightarrow \Hom_\T(M,Y[i])\rightarrow\Hom_\T(K,Y[i+1]).$$ Therefore $\Hom_\T(M,Y[i])=0$ for $i>\alpha$ since $\pd_{\Y}\,(K)\leq \alpha.$ So we get that $\pd_{\Y}\,(\X^\vee) \leq\pd_{\Y}\,(\X).$
\

Finally, it is easy to see that (b) is a consequence of (a).
\end{dem}
\vspace{.2cm}

The following two lemmas resembles the so called ``shifting argument'' that is usually used for syzygies and cosyzygies in the $\Ext^n$ functor.

\begin{lem}\label{corricores} Let $\X$ and $\Y$ be classes of objects in $\T$ such that $\id_{\X}(\Y)=0.$ Then, for any $X\in\X,$ $k>0$ and $K_n\in\Y*\Y[1]*\cdots*\Y[n-1]*K_0[n],$ there is an isomorphism of abelian groups
 $$\Hom_\T(X,K_0[k+n])\simeq\Hom_\T(X,K_n[k]).$$
\end{lem}
\begin{dem} Let $X\in\X,$ $k>0$ and $K_n\in\Y*\Y[1]*\cdots*\Y[n-1]*K_0[n].$ By \ref{restri} (a), we have distinguished triangles $\eta_i:\quad K_i\rightarrow Y_i\rightarrow K_{i+1}\rightarrow K_i[1]$ with $Y_i\in\Y,\;0\leq i\leq n-1.$ Applying the functor $\Hom_\T(X[-k],-)$ to $\eta_i,$ we get the exact sequence of abelian groups
 $$(X[-k],Y_i)\rightarrow(X[-k],K_{i+1})\rightarrow(X[-k],K_i[1])\rightarrow(X[-k],Y_i[1]),$$
where $(-,-):=\Hom_\T(-,-)$ for simplicity. Since $\id_{\X}(\Y)=0,$ it follows that $\Hom_\T(X[-k],K_{i+1})\simeq\Hom_\T(X[-k],K_i[1]).$ Therefore, by the preceding isomorphism, we have\\ $\qquad\Hom_\T(X,K_{n}[k])\simeq\Hom_\T(X,K_{n-1}[k+1])\simeq\cdots\simeq\Hom_\T(X,K_0[k+n]).$
\end{dem}

\begin{lem}\label{corrires} Let $\X$ and $\Y$ be classes of objects in $\T$ such that $\pd_{\X}(\Y)=0.$ Then, for any $X\in\X,$ $k>0$ and $K_n\in K_0[-n]*\Y[-n+1]*\cdots*\Y[-1]*\Y,$ there is an isomorphism of abelian groups
 $$\Hom_\T(K_0,X[k+n])\simeq\Hom_T(K_n,X[k]).$$
\end{lem}
\begin{dem} The proof is similar to the one given in \ref{corricores} by using \ref{restri} (b).
\end{dem}

\section{relative weak-cogenerators and relative injectives}

In this section, we focus our attention on pairs $(\X,\omega)$ of classes of objects in $\T.$  We study the relationship between weak-cogenerators in $\X$ and coresolutions. Also, we give a characterization of some special subcategories of $\T.$

\begin{defi} Let $(\X,\omega)$ be a pair of classes of objects in $\T.$ We say that
\begin{itemize}
 \item[(a)] $\omega$ is a {\bf{weak-cogenerator}} in $\X,$ if $\omega\subseteq\X\subseteq\X[-1]*\omega;$
 \item[(b)] $\omega$ is a {\bf{weak-generator}} in $\X,$ if $\omega\subseteq\X\subseteq \omega * \X[1];$
 \item[(c)] $\omega$ is {\bf{$\X$-injective}} if $\id_{\X}(\omega)=0;$ and dually,  $\omega$ is {\bf{$\X$-projective}}  if $\pd_{\X}(\omega)=0.$
\end{itemize}
\end{defi}

The following result say us that an $\X$-injective weak-cogenerator, closed under direct summands, is unique (in case it exists).

\begin{pro}\label{coinyec0} Let $(\X,\omega)$ be a pair of classes of objects in $\T$ such that $\omega$ is $\X$-injective.
Then, the following statements hold.
 \begin{enumerate}
  \item[(a)] $\omega^\wedge$ is $\X$-injective.
  \item[(b)] If $\omega$ is a weak-cogenerator in $\X,$ and $\omega$ is closed
under direct summands in $\T,$ then $$\omega=\X\cap{}_{\X}\U^{\perp}[-1]=\X\cap\omega^\wedge.$$
\end{enumerate}
\end{pro}
\begin{dem} (a) It follows from the dual result of \ref{pdim} (a).
\

(b) Let $\omega\subseteq\X\subseteq\X[-1]*\omega$ and $\omega$ be closed
under direct summands in $\T.$
\

 We start by proving the first equality. Let $X\in \X\cap{}_{\X}\U^{\perp}[-1]$. Since $\X\subseteq\X[-1]*\omega,$ there is a distinguished triangle
$$\eta:\quad X\rightarrow W\rightarrow X'\stackrel{f}{\rightarrow} X[1]\text{ in $\T$ with $X'\in\X$ and $W\in\omega.$}$$
Moreover $X\in{}_{\X}\U^{\perp}[-1]$ implies that $\Hom_\T(-,X[1])|_\X=0$ (see \ref{susp-perp} (b)). Hence $\eta$ splits and so $X\in\omega;$ proving that $\X\cap{}_{\X}\U^{\perp}[-1]\subseteq\omega.$ The other inclusion follows from \ref{compdim} (a2) since $\omega\subseteq\X$ and $\id_{\X}(\omega)=0.$\\
On the other hand,  it is easy to see that $\omega\subseteq \X\cap\omega^\wedge$  and
since $\id_{\X}(\omega^\wedge)=0,$ it follows from \ref{compdim} (a2) that $\X\cap\omega^\wedge\subseteq \X\cap{}_{\X}\U^{\perp}[-1];$ proving (b).
\end{dem}

\begin{pro}\label{coresIdfin} Let $(\X,\omega)$ be a pair of classes of objects in $\T,$ and $\omega$ be closed under direct summands in $\T.$ If $\omega$ is an $\X$-injective weak-cogenerator in $\X,$ then
$$\X\cap\omega^\vee=\{X\in\X\,:\,\id_{\X}(X)<\infty\}.$$
\end{pro}
\begin{dem} Let $M\in\X\cap\omega^\vee.$ We assert that $\id_\X(M)\leq d<\infty$ where $d:=\coresdim_\omega(M).$ Indeed, from \ref{rkepsilon} (a), dual version, and \ref{restri} (a), there is some $W_d\in\omega*\omega[1]*\cdots*\omega[d-1]*M[d]$ with $W_d\in\omega.$ So, by \ref{corricores} we get an isomorphism $\Hom_\T(X,M[k+d])\simeq\Hom_\T(X,W_d[k])$ for any $k>0$
and $X\in\X;$ and using that $\id_{\X}(\omega)=0,$ it follows that $\Hom_\T(X,M[k+d])=0$ for any $k>0,$ proving that $\id_\X(M)\leq d.$
\

Let $N\in\X$ be such that $n:=\id_\X(N)<\infty.$ Using that $\X\subseteq\X[-1]*\omega,$ we can construct a family $\{K_i\rightarrow W_i\rightarrow K_{i+1}\rightarrow K_i[1]\;:\; W_i\in\omega,\;0\leq i\leq n-1\}$ of distinguished triangles in $\T$ where $K_0:=N$ and $K_i\in\X,$ $\forall i\;0\leq i\leq n.$ Thus, by \ref{restri} (a), it follows that $K_n\in\omega*\omega[1]*\cdots*\omega[n-1]*N[n];$ and so by \ref{corricores} we get that
$\Hom_T(X,K_n[k])\simeq \Hom_\T(X,N[k+n]),$ $\forall X\in\X,$ $\forall k>0.$ But $  \Hom_\T(X,N[k+n])=0,$ $\forall X\in\X,$ $\forall k>0$ because $\id_{\X}(N)=n.$  Therefore $\id_\X(K_n)=0$ and then $K_n\in\omega$ (see \ref{compdim} and \ref{coinyec0} (b)); proving that $N\in\X\cap\omega^\vee.$
\end{dem}
\

Now, we are in condition to prove the following result. In the statement, we use the notions of precovering and preenveloping classes (see Section 1).
\begin{teo}\label{specialtrian} Let $(\X,\omega)$ be a pair of classes of objects in $\T,$ $\X$ be closed under extensions and $\omega$ be a weak-cogenerator in $\X.$ Then, the following statements hold.
\begin{enumerate}
 \item[(a)] For all $C\in\X^\wedge$ there exist two distinguished triangles in $\T:$\\
  $\begin{CD} C[-1] @>>> Y_C @>>> X_C @>\varphi_C>> C \text{ with }Y_C\in\omega^\wedge\text{ and } X_C\in\X,\end{CD}$
\noindent  $\begin{CD} C @>\varphi^C>> Y^C @>>> X^C @>>> C[1] \text{ with }Y^C\in\omega^\wedge\text{ and } X^C\in\X. \end{CD}$
 \item[(b)] If $\omega$ is $\X$-injective, then
\vspace{.2cm}
  \begin{itemize}
   \item[(b1)] $Y_C[1]\in\X^\perp$ and $\varphi_C$ is an  $\X$-precover of $C,$
\vspace{.2cm}
   \item[(b2)] $X^C[-1]\in {}^\perp(\omega^\wedge)$ and $\varphi^C$ is a $\omega^\wedge$-preenvelope of $C.$
  \end{itemize}
\end{enumerate}
\end{teo}
\begin{dem} (a) Let $C\in\X^\wedge.$  We prove the existence of the triangles in (a) by induction on $n:=\resdim_\X(C).$ If $n=0,$ we have that $C\in\X$ and then we can consider $C[-1] \rightarrow 0\rightarrow C\stackrel{1_C}{\rightarrow} C$ as the first triangle; the second one can be obtained from the fact that $\X\subseteq\X[-1]*\omega.$
\

Assume that $n>0.$ Then, we have a distinguished triangle $C[-1]\rightarrow K_1\rightarrow X_0\rightarrow C$ in $\T$ with $X_0\in\X$ and $\resdim_\X(K_1)=n-1.$ Hence, by induction, there is a distinguished triangle
 $ K_1 \to Y^{K_1} \to X^{K_1} \to K_1[1] $ in $\T$ with $Y^{K_1}\in\omega^\wedge$ and $X^{K_1}\in\X.$ By the co-base change procedure applied to  the above triangles, there exists a commutative diagram
$$\begin{CD}
@. X^{K_1}[-1] @= X^{K_1}[-1] @.\\
@. @VVV @VVV @.\\
C[-1] @>>> K_1 @>>> X_0 @>>> C\\
@| @VVV @VVV @|\\
C[-1] @>>> Y^{K_1} @>>> U @>>> C\\
@. @VVV @VVV @.\\
@. X^{K_1} @= X^{K_1} @.\\
{}
  \end{CD}$$
where the rows and columns are distinguished triangles in $\T.$ Since $X_0,X^{K_1}\in\X$ it follows that $U\in\X.$ By taking $X_C:=U$ and $Y_C:=Y^{K_1},$ we get the first triangle in (a). On the other hand, since $U\in\X$ and $\X\subseteq\X[-1]*\omega,$ there exists  a distinguished triangle $X^C[-1]\rightarrow U\rightarrow W\rightarrow X^C$ in $\T$  with $X^C\in \X$  and $W\in\omega.$ Again, by the co-base change procedure, there exists a commutative diagram
$$\begin{CD}
@. Y^{K_1} @= Y^{K_1} @.\\
@. @VVV @VVV @.\\
X^C[-1] @>>> U @>>> W @>>> X^C\\
@| @VVV @VVV @|\\
X^C[-1] @>>> C @>>> Y^C @>>> X^C\\
@. @VVV @VVV @.\\
@. Y^{K_1}[1] @= Y^{K_1}[1] @.\\
{}
  \end{CD}$$
where the rows and columns are distinguished triangles in $\T.$ By the second column, in the diagram above,
it follows that $Y^C\in\omega^\wedge.$ Hence the second row in the preceding diagram is the desired triangle.
\

(b2) Consider the triangle $X^C[-1]\stackrel{g}{\rightarrow}
C\stackrel{\varphi^C}{\to} Y^C\to X^C$ with $Y^C\in\omega^\wedge$
and $X^C\in\X.$ Since $\id_{\X}(\omega)=0$ we have by \ref{coinyec0}
that $\id_{\X}(\omega^\wedge)=0.$ Thus
$\Hom_\T(X[-1],-)|_{\omega^\wedge}=0$ for any $X\in\X;$ and so
$X^C[-1]\in {}^\perp(\omega^\wedge).$ Let $f:C\to Y$ be a morphism
in $\T$ with $Y\in\omega^\wedge.$ Since $\Hom_\T(X^C[-1],Y)=0,$ we
have that $fg=0$ and hence $f$ factors through $\varphi^C;$ proving
that $\varphi^C$ is a $\omega^\wedge$-preenvelope of $C.$ \

(b1) It is similar to the proof of (b2).
\end{dem}
\vspace{.2cm}

The following result provides a characterization of the category $\X^\wedge,$ and will be applied in \cite{MSSS2} to deal with co-$t$-structures.

\begin{cor}\label{igualresdim} Let $(\X,\omega)$ be a pair of classes of objects in $\T$ such that $\X$ is closed under extensions and $\omega$ is a weak-cogenerator in $\X.$ Then, the following statements hold.
\begin{enumerate}
  \item[(a)] If $0\in \omega$ then $\X^\wedge=\X*\omega^\wedge=\X*\omega^\wedge [1].$
\vspace{.2cm}
  \item[(b)]  If $\X[-1]\subseteq\X $ then $\X^\wedge=\X*\omega^\wedge=\X*\omega^\wedge [1]=\X[-1]*\omega^\wedge .$
 \end{enumerate}
\end{cor}
\begin{dem} We assert that $\X*\omega^\wedge \subseteq \X^\wedge.$ Indeed, since $\omega\subseteq\X$ it follows from \ref{rkepsilon} (b) that $\varepsilon^\wedge_n(\omega)\subseteq \varepsilon^\wedge_n(\X),$ giving us that $\omega^\wedge \subseteq \X^\wedge.$ Hence
$\X*\omega^\wedge \subseteq \X*\X^\wedge$ and then $\X*\omega^\wedge \subseteq \X^\wedge $  by \ref{stardim} (b).
\

(a) Let $0\in \omega.$ By \ref{specialtrian} (a) we have that $\X^\wedge\subseteq\X*\omega^\wedge[1],$ and therefore, by \ref{rkepsilon} (d)  we get  $\X^\wedge\subseteq \X*\omega^\wedge[1] \subseteq \X*\omega^\wedge.$ But $\X*\omega^\wedge   \subseteq \X*\X^\wedge   \subseteq   \X^\wedge  $  by  \ref{stardim} (b),   and then
  $\X^\wedge=\X*\omega^\wedge=\X*\omega^\wedge [1].$
\

(b) Let $\X[-1]\subseteq\X .$ By \ref{specialtrian} (a) and the assertion above,  we have
$\X^\wedge\subseteq \X [-1]*\omega^\wedge \subseteq \X*\omega^\wedge  \subseteq \X^\wedge.$ On the other hand, from \ref{specialtrian} (a), it follows that $\X^\wedge\subseteq\X*\omega^\wedge[1].$ So, to prove (b), it is enough to see that $\X*\omega^\wedge[1]\subseteq \X^\wedge.$ Let $C\in\X*\omega^\wedge[1].$ Then there is a distinguished triangle $Y \to X \to C \to Y[1]$ in $\T$ with $X\in\X$ and $Y\in\omega^\wedge.$ Hence it follows that $C\in\X^\wedge$ since $\omega^\wedge \subseteq \X^\wedge;$ proving (b).
\end{dem}
\vspace{.2cm}

We are now in position to prove that if $\omega$ is an $\X$-injective weak-cogenerator in a suitable class $\X,$ then the $\omega^\wedge$-projective dimension coincides with the $\X$-resolution dimension for every object of the thick subcategory of $\T$ generated by $\X.$

\begin{teo}\label{pd-resd} Let $(\X,\omega)$ be a pair of classes of objects in $\T$ which are closed under direct summands in $\T.$ If $\X$ is closed under extensions and $\omega$ is an $\X$-injective weak-cogenerator in $\X,$ then
 $$\pd_{\omega^\wedge}(C)=\pd_\omega(C)=\resdim_{\X}(C),\quad\forall C\in\X^\wedge.$$
\end{teo}
\begin{dem}  Let $C\in\X^\wedge.$ By \ref{compdim} (b) and the dual of \ref{pdim} (a), it follows that
 $\pd_\omega(C)=\id_{\{C\}}(\omega)=\id_{\{C\}}(\omega^\wedge)=\pd_{\omega^\wedge}(C).$
To prove the last equality, we proceed by induction on $n:=\resdim_\X(C).$ To start with, we have $\pd_\omega(\X)=\id_\X(\omega)=0.$ If $n=0$ then $C\in\X$ and so $\pd_\omega(C)=0=\resdim_{\X}(C).$
\

Let $n=1.$ Then, we have a distinguished triangle $X_1\to X_0\to C\to X_1[1]$ in $\T$ with $X_i\in\X.$ By
\ref{specialtrian} (a), there is a distinguished triangle $Y_C \to X_C \stackrel{\varphi_C}{\to} C \to Y_C[1]$ in $\T$ with $Y_C\in\omega^\wedge$ and $X_C\in\X.$ By the base change procedure, there exists a commutative diagram\\
$$\begin{CD}
@. Y_C @= Y_C @.\\
@. @VVV @VVV @.\\
X_1 @>>> E @>>> X_C @>>> X_1[1]\\
@| @VVV @V\varphi_C VV @|\\
X_1 @>>> X_0 @>\alpha >> C @>>> X_1[1]\\
@. @VVV @V\beta VV @.\\
@. Y_C[1] @= Y_C[1]\,, @.\,
  \end{CD}$$
\noindent
${}$
\vspace{.2cm}

where the rows and columns are distinguished triangles in $\T.$ Since $X_1,X_C\in\X$ it follows that $E\in\X.$ On the other hand, since $\Hom_\T(X,Y[1])=0$ for any $X\in\X$ and $Y\in\omega^\wedge$ (see \ref{coinyec0} (a)), we get that $\beta\alpha=0$ and then the triangle $Y_C\to E\to X_0\to Y_C[1]$ splits getting us that $Y_C\in\X\cap\omega^\wedge=\omega$ (see \ref{coinyec0}). Using that $\pd_\omega(\X)=0$ and \ref{pdyresdim} (a), we have that $\pd_\omega(C)\leq\resdim_\X(C)=1.$ We assert that $\pd_\omega(C)>0.$ Indeed, suppose that $\pd_\omega(C)=0;$ and then $\Hom_\T(C,W[1])=0$ for any $W\in\omega.$ Since $Y_C\in\omega$ we get that $\beta=0$ and hence the triangle $Y_C\to X_C\to C\to Y_C[1]$ splits. Therefore $C\in\X$ contradicting that $\resdim_\X(C)=1;$ proving that $\pd_\omega(C)=1=\resdim_\X(C).$
\

Let $n\geq 2.$ From \ref{pdyresdim} (a), we have that $\pd_\omega(C)\leq\resdim_\X(C)=n$ since $\pd_\omega(\X)=0.$ Then, it is enough to prove that $\Hom_\T(C[-n],-)|_\omega\neq 0.$ Consider a distinguished triangle $K_1\to X_0\to C\to K_1[1]$ in $\T$ with $X_0\in\X$ and $\resdim_\X(K_1)=n-1=\pd_\omega(K_1).$ Applying the functor $\Hom_\T(-,W[n]),$ with $W\in\omega,$ to the triangle $C\to K_1[1]\to X_0[1]\to C[1]$ we get the exact sequence of abelian groups
$$\Hom_\T(X_0[1],W[n])\to\Hom_\T(K_1[1],W[n])\to\Hom_\T(C,W[n]).$$ Suppose that $\Hom_\T(C[-n],-)|_\omega=0.$ Then $\Hom_\T(K_1[1],W[n])=0$ since $\id_{\X}(\omega)=0$ and $n\geq 2;$ contradicting that $\pd_\omega(K_1)=n-1.$
\end{dem}

\begin{lem}\label{relid} Let $\X$ be a class of objects in $\T$ and $A\rightarrow B\rightarrow C\rightarrow A[1]$ a distinguished triangle in $\T.$ Then
\begin{enumerate}
 \item[(a)] $\id_{\X}(B)\leq\maxi\,\{\id_{\X}(A),\id_{\X}(C)\};$
\vspace{.2cm}
 \item[(b)] $\id_{\X}(A)\leq\maxi\,\{\id_{\X}(B),\id_{\X}(C)+ 1 \};$
\vspace{.2cm}
 \item[(c)] $\id_{\X}(C)\leq\maxi\,\{\id_{\X}(B),\id_{\X}(A)-1 \}.$
\end{enumerate}
\end{lem}
\begin{dem} It is straightforward.
\end{dem}
\vspace{.2cm}

The following result gives a relationship between the relative injective dimensions, attached to
the pair $(\X,\omega),$ and the $\omega$-coresolution dimension. Such a result will be applied in 
\cite{MSSS2} to deal with co-$t$-structures. Observe that \ref{idXcores} is not the dual version 
of \ref{pd-resd}.

\begin{pro}\label{idXcores} Let $(\X,\omega)$ be a pair of classes of objects in $\T$ such that $\omega\subseteq{}_\X\U.$ If $\omega$ is closed under direct summands and $\X$-injective, then
 $$\id_\omega(C)=\id_\X(C)=\coresdim_\omega(C),\quad\forall\,C\in{}_\X\U\cap\omega^\vee.$$
\end{pro}
\begin{dem} Assume that $\omega$ is closed under direct summands and $\id_{\X}(\omega)=0.$ Let $C\in{}_\X\U\cap\omega^\vee$ and $n:=\coresdim_\omega(C).$ By the dual of \ref{pdyresdim} (b), it follows $$(*)\quad\alpha:=\id_\omega(C)\leq\id_{\X}(C)=\coresdim_\omega(C)=n.$$ Moreover, since $C\in\omega^\vee$ there is a distinguished triangle
$(\eta)\;:\;C\rightarrow W_0\rightarrow K_1\rightarrow C[1]$ in $\T$ with $W_0\in\omega$ and $\coresdim_\omega(K_1)=n-1.$ Furthermore, from \ref{rkepsilon} (a) we get that $K_1\in{}_\X\U$ since ${}_\X\U$ is closed under cocones and $\omega\subseteq{}_\X\U.$ Now, we prove the result by induction on $\alpha.$
\

Let $\alpha=0.$ We assert that $C\in\omega$ (note that if this is true, then the result follows). We proceed by induction on $n.$ If $n=0$ it is clear that $C\in\omega.$ So we may assume that $n>0,$ and then, applying \ref{relid} to $(\eta)$ it follows that $\id_\omega(K_1)=0.$ Hence by induction we get that $K_1\in\omega,$ and so $\Hom_\T(K_1,C[1])=0$ since $\id_\omega(C)=0.$ Therefore the triangle $(\eta)$ splits and then $C\in\omega;$ proving the assertion.
\

Assume that $\alpha>0.$ Applying \ref{relid} to $(\eta),$ we get that $\id_\omega(K_1)\leq\alpha-1.$ Thus, by induction, it follows that $\id_\omega(K_1)=\id_\X(K_1)=\coresdim_\omega(K_1)=n-1.$ In particular, we obtain that $n-1\leq\alpha-1$ and hence by $(*)$ the result follows.
\end{dem}
\vspace{.2cm}

The following result provides a characterization of $\omega^\wedge$ in terms of $\X.$ From this, we get some nice relationship between other subcategories.

\begin{pro}\label{resomega} Let $(\X,\omega)$ be a pair of classes of objects in $\T$ such that $\omega$ is closed under direct summands in $\T,$ $\X$ is closed under extensions and $\omega$ is an $\X$-injective weak-cogenerator in $\X.$ Then, the following statements hold.
 \begin{itemize}
  \item[(a)] ${}_\X\U^\perp[-1]\cap\X^\wedge=\omega^\wedge.$
\vspace{.2cm}
  \item[(b)] If $\X[-1]\subseteq\X$ then $\U_\omega=\omega^\wedge=\X^\perp[-1]\cap\X^\wedge.$
 \end{itemize}
\end{pro}
\begin{dem} (a) Let $C\in{}_\X\U^\perp[-1]\cap\X^\wedge.$ In particular, from \ref{specialtrian} (a), there exists a distinguished triangle $Y_C\to X_C\to C\to Y_C[1]$ in $\T$ with $Y_C\in\omega^\wedge$ and $X_C\in\X.$ We assert that $\id_\X(X_C)=0.$ Indeed, it follows from \ref{relid} (a) since $\id_\X(C)=0=\id_\X(Y_C)$ (see \ref{compdim} and \ref{coinyec0} (a)). Therefore, $X_C\in\X\cap{}_\X\U^\perp[-1]$ and by \ref{coinyec0} (b), we get that $X_C\in\omega$ proving that $C\in\omega^\wedge.$ On the other hand, since $\id_\X(\omega^\wedge)=0,$  we have from \ref{compdim} that $\omega^\wedge\subseteq{}_\X\U^\perp[-1]\cap\X^\wedge.$
\

(b) Assume that $\X[-1]\subseteq\X.$ Hence, by \ref{xx=x} (b), we have that $\X$ is a cosuspended subcategory of $\T.$ Therefore, from (a), it follows that $\omega^\wedge=\X^\perp[-1]\cap\X^\wedge.$ Furthermore, since $\X^\perp[-1]$ is suspended and $\X^\wedge$ is triangulated (see \ref{resBuan}), we conclude that $\omega^\wedge$ is a suspended subcategory of $\T;$ and so $\U_\omega\subseteq\omega^\wedge.$ Finally, the equality $\U_\omega=\omega^\wedge$ follows from \ref{varios} (3).
\end{dem}

\begin{teo}\label{thicksub} Let $(\X,\omega)$ be a pair of classes of objects in $\T$ which are closed under direct summands, $\X$  be cosuspended and $\omega$ be an $\X$-injective weak-cogenerator in $\X.$ Then,
 $$\varepsilon^\wedge_n(\X)=\X[n]=\X^\wedge\cap{}^\perp\U_\omega[n+1]=\X^\wedge\cap{}^\perp(\omega^\wedge)[n+1],\quad\forall n\geq 0.$$
\end{teo}
\begin{dem} From \ref{resBuan}, we have that $\varepsilon^\wedge_n(\X)=\X[n]$ and $\X^\wedge=\cup_{n\geq 0}\,\X[n].$ On the other hand, by \ref{compdim} and \ref{pd-resd}, it follows that $$\X^\wedge\cap{}^\perp\U_{\omega^\wedge}[n+1]=\X^\wedge\cap{}^\perp\U_\omega[n+1]=\X[n]\cap\X^\wedge=\X[n].$$ Finally, since $\omega^\wedge$ is a suspended subcategory of $\T$ (see \ref{resomega} (b)), we have that ${}^\perp\U_{\omega^\wedge}={}^\perp(\omega^\wedge);$ proving the result.
\end{dem}

\vspace{.2cm}

The previous results can be seen under the light of the so called torsion theories, in the sense of 
Iyama-Yoshino. Such torsion theories have been extensively studied in relation to the cluster theory (see \cite{IY}).

\begin{defi}\cite[Definition 2.2]{IY} A pair $(\X,\Y)$ of subcategories of $\T$ is called a torsion theory in $\T,$ if the following conditions hold.
\begin{itemize}
 \item[(a)] $\X$ and $\Y$ are closed under direct summands in $\T.$
 \item[(b)] $\Hom_\T(\X,\Y)=0.$
 \item[(c)] $\T=\X*\Y.$
\end{itemize} 
\end{defi}

\begin{cor}\label{torsion} Let $(\X,\omega)$ be a pair of classes of objects in $\T,$ which are closed under 
direct summands in $\T,$ and such that $\X$ is 
cosuspended and $\omega$ is an $\X$-injective weak-cogenerator in $\X.$ Then, the pair 
$(\X,\omega^\wedge [1])$ is a torsion theory in the thick triangulated subcategory $\X^\wedge$ of
 $\T.$
\end{cor}
\begin{dem} Since $\X$ is cosuspended and closed under direct summands, we know from 
\ref{resBuan} (c) that $\X^\wedge$ is a thick triangulated subcategory of $\T.$ On the other 
hand, from \ref{igualresdim}, it follows that $\X^\wedge=\X*\omega^\wedge [1];$ and furthermore, 
$\Hom_\T(\X,\omega^\wedge [1])=0$ since $\omega$ is  $\X$-injective. On the other hand, by 
\ref{resomega}, we get that  $\omega^\wedge$ is a subcategory  closed under direct summands in 
$\T.$ 
\end{dem}
\vspace{.2cm}

\begin{defi} For a given class $\Y$ of objects in $\T,$ we set $\Y^\sim:=(\Y^\wedge)^\vee.$
\end{defi}

\begin{lem}\label{Xsim} Let $\X$ be a class of objects in $\T.$ Then, the following statements hold.
 \begin{itemize}
  \item[(a)] If $\X^\wedge$ is closed under cocones then $\omega^\sim\subseteq\X^\wedge$ for any $\omega\subseteq\X.$
\vspace{.2cm}
  \item[(b)] $\X^\wedge$ is closed under cocones if and only if $\X^\wedge=\X^\sim.$
\vspace{.2cm}
  \item[(c)] If $\X^\wedge=\X^\sim$ then $\X^\wedge[-1]\subseteq\X^\wedge.$
 \end{itemize}
\end{lem}
\begin{dem} (a) Let $\omega\subseteq\X$ and assume that $\X^\wedge$ is closed under cocones. Hence $\omega^\wedge\subseteq\X^\wedge$ and so by \ref{varios} (2), we conclude that $\omega^\sim\subseteq\X^\wedge.$
\

(b) Assume that $\X^\wedge$ is closed under cocones. It is clear that $\X^\wedge\subseteq\X^\sim.$ On the other hand, by (a) it follows that  $\X^\sim\subseteq\X^\wedge.$
\

Suppose that $\X^\wedge=\X^\sim.$ Let $A\to B\to C\to A[1]$ be a distinguished triangle in $\T$ with $B,C$ in $\X^\wedge.$ Then $A\in\X^\sim=\X^\wedge$ and so $\X^\wedge$ is closed under cocones.
\

(c) Let $\X^\wedge=\X^\sim$ and consider $X\in\X^\wedge.$ Since, we have the distinguished triangle $X[-1]\to 0\to X\stackrel{1_X}{\to} X$ and $0,X \in \X^{\wedge},$ it follows from (b) that $X[-1]\in\X^\wedge;$ proving the lemma.
\end{dem}

\begin{cor}\label{tilde-eq} Let $(\X,\omega)$ be a pair of classes of objects in $\T.$ If $\X$ is cosuspended and $\omega\subseteq\X,$ then $\omega^\sim\subseteq\X^\wedge=\X^\sim.$
\end{cor}
\begin{dem} It follows from \ref{Xsim} and the fact that $\X^\wedge$ is triangulated (see \ref{resBuan}).
\end{dem}

In case $\omega$ is an $\X$-injective weak-cogenerator in a cosuspended subcategory $\X$ of $\T,$ both closed under direct summands, the
thick subcategory $\overline{\Delta}_\T(\omega)$ of $\T$ can be characterized as follows.

\begin{teo}\label{tilde-id} Let $(\X,\omega)$ be a pair of classes of objects in $\T,$  $\X$ be cosuspended and $\omega$ be closed under direct summands in $\T.$ If $\omega$ is an $\X$-injective weak-cogenerator in $\X,$ the following statements hold.
 \begin{itemize}
  \item[(a)] $\omega^\sim=\{C\in\X^\wedge : \id_{\X}(C)<\infty\}=\X^\wedge\cap(\X^\perp[-1])^\vee.$
\vspace{.2cm}
  \item[(b)] $\omega^\sim$ is the smallest triangulated subcategory of $\X^\wedge$ containing $\omega,$ that is
$\omega^\sim=\Delta_{\X^\wedge}(\omega).$
\vspace{.2cm}
  \item[(c)] If $\X$ is closed under direct summands in $\T,$ then $$\overline{\Delta}_\T(\omega)=\omega^\sim=\overline{\Delta}_\T(\X)\cap(\X^\perp[-1])^\vee.$$
 \end{itemize}
\end{teo}
\begin{dem} Assume that $\omega\subseteq\X\subseteq \X[-1]*\omega$ and $\id_{\X}(\omega)=0.$ Let $\Y:=\{C\in\X^\wedge : \id_{\X}(C)<\infty\}.$ We start by proving that $\omega^\sim\subseteq\Y.$ By \ref{tilde-eq}, we know that $\omega^\sim\subseteq\X^\wedge.$ On the other hand, since $\id_{\X}(\omega^\wedge)=0$ (see \ref{coinyec0}(a)), we can apply the dual of $\ref{pdyresdim} (a),$ and then $\id_\X(C)\leq\coresdim_{\omega^\wedge}(C)<\infty$ for any $C\in\omega^\sim;$ proving that $\omega^\sim\subseteq\Y.$\\
Let $C\in\Y.$ By \ref{specialtrian} (a), there is a distinguished triangle $C\to Y^C\to X^C\to C[1]$ in $\T$ with $Y^C\in\omega^\wedge$ and $X^C\in\X.$ Hence, from \ref{relid} (c) we get that $\id_{\X}(X^C)<\infty$ and then, by \ref{coresIdfin} $X^C\in\omega^\vee \subseteq \omega^\sim;$ proving that $C\in\omega^\sim.$ Hence $\Y\subseteq\omega^\sim.$ In order to get the second equality in (a), we use \ref{compdim} and the fact that $\X={}_\X\U$ to obtain
 $$\{C\in\X^\wedge : \id_{\X}(C)<\infty\}=\X^\wedge\cap(\cup_{n\geq 0}\,\X^\perp[-n-1]).$$ On the other hand, since $\X^\perp[-1]$ is suspended, then by the dual of \ref{resBuan}, it follows that $(\X^\perp[-1])^\vee=\cup_{n\geq 0}\,\X^\perp[-n-1]$  and also that $(\X^\perp[-1])^\vee$ is a thick subcategory of $\T.$ In particular, by \ref{resBuan}, we get (b). Finally, (c) follows from (a) and \ref{resBuan}.
\end{dem}
\vspace{.2cm}

The following result will be applied in \cite{MSSS2} to deal with co-$t$-structures. Also, it will be applied in Section 5, to get some connection
with relative Rouquier's dimension.

\begin{pro}\label{idwX} Let $(\X,\omega)$ be a pair of classes of objects in $\T,$ $\X$ cosuspended and $\omega$ closed under direct summands in $\T.$ If $\omega$ is an $\X$-injective weak-cogenerator in $\X,$ then
 \begin{itemize}
 \item[(a)] $\id_\omega(C)=\id_\X(C)<\infty,\quad\forall C\in\omega^\sim;$
\vspace{.2cm}
 \item[(b)] $\omega^\sim\cap{}_\omega\U^\perp[-n-1]=\omega^\sim\cap\X^\perp[-n-1],\quad \forall n\geq 0.$
 \end{itemize}
\end{pro}
\begin{dem} (a) By \ref{resBuan} and \ref{tilde-id}, we know that $\X^\wedge$ and $\omega^\sim$ are triangulated subcategories of $\T.$ Furthermore, from \ref{tilde-eq} it follows that $\omega^\sim\subseteq\X^\wedge.$ Let $C\in\omega^\sim.$ It is enough to prove that $\id_\X(C)\leq\id_\omega(C).$ In order to do that, we will use induction on $n:=\id_\omega(C).$
\

Since $C\in\X^\wedge,$ we have from \ref{specialtrian} the existence of a distinguished triangle $(\eta)\;:\;C\to Y^C\to X^C\to C[1]$ in $\T$ with $Y^C\in\omega^\wedge\subseteq\omega^\sim$ and $X^C\in\X.$ We assert that $X^C\in\X\cap\omega^\vee.$ Indeed, using that $\omega^\sim$ is triangulated we conclude that $X^C\in\X\cap\omega^\sim$ and hence $\id_\X(X^C)$ is finite (see \ref{tilde-id} (a)). Thus $X^C\in\X\cap\omega^\vee$ by \ref{coresIdfin}; proving the assertion.
\

Let $n=0.$ Then $\id_\omega(X^C)=0$ since $\id_\omega(Y^C)=0$ (see \ref{coinyec0} and \ref{relid}). On the other hand, \ref{idXcores} gives the equalities $\coresdim_\omega(X^C)=\id_\omega(X^C)=0.$ Hence $X^C\in\omega$ and since $\id_\omega(C)=0,$ it follows that $\Hom_\T(X^C,C[1])=0.$ Therefore, the triangle $(\eta)$ splits giving us that $C$ is a direct summand of $Y^C,$ and hence $\id_\X(C)\leq\id_\X(Y^C)\leq\id_\X(\omega^\wedge)=0.$
\

Assume that $n>0.$ Since $\id_\X(Y^C)=0=\id_\omega(Y^C),$ it follows from \ref{relid} that $\id_\omega(X^C)\leq n-1.$ Hence, by induction $\id_\X(X^C)\leq\id_\omega(X^C)\leq n-1.$ Therefore, applying again \ref{relid} to the triangle $(\eta),$ we get that $\id_\X(C)\leq n=\id_\omega(C);$ proving the result.
\

(b) By \ref{compdim}, the item (a) and the fact that ${}_\X\U=\X$ the result follows.
\end{dem}

\section{Some connection with Rouquier's dimension}

In this section, we introduce some kind of ``relative Rouquier's dimension''; and relate it with the Rouquier's dimension and the other relative 
dimensions developed in this paper.
\

Let $\X$ and $\Y$ be classes of objects in a triangulated category $\T.$
Consider  the smallest subcategory $\left\langle \X\right\rangle$ of $\T$ containing $\X,$ closed under shifts,  finite direct sums and direct summands, that is,  $\left\langle \X\right\rangle:=\add\,(\cup_{i\in\Enteros}\,\X[i]).$ Let  
$\X\diamondsuit\Y:=\left\langle\X*\Y\right\rangle.$ Following R. Rouquier, we inductively define 
$\la \X\ra_0:=0$ and $\la \X\ra_n:=\la \X\ra_{n-1}\diamondsuit\la \X\ra$ for $n\geq 1.$ 
The objects of  $\la \X\ra_{n}$ are the direct summands of the objects obtained by taking a n-fold
extension of finite direct sums of shifts of objects of $\X.$

The following dimension for triangulated categories was introduced by R. Rouquier and has been extensively studied in \cite{Ro}.

\begin{defi}\cite[Definition 3.2]{Ro} Let $\T$ be a triangulated category. The dimension of $\T$ is $$\dim\,(\T):=\mini\{n\in\N\;\;|\;\text{ there exists }X\in\T\text{ such that }
\la X\ra_{n+1}=\T\}.$$
\end{defi}

\begin{lem}\label{conexion} Let $\X$ be a class of objects in a triangulated category $\T.$ Then, for any $n\in\N,$ the following statements hold.
\begin{itemize}
 \item[(a)] $\la \X\ra_{n+1}=\la\varepsilon_n^\wedge(\la \X\ra)\ra=\la\varepsilon_n^\vee(\la \X\ra)\ra.$
 \item[(b)] $\la \X\ra_{n}\subseteq\la \X\ra_{n+1}.$
\end{itemize}
\end{lem}
\begin{dem} (a) By induction over $n,$ it can be seen that 
$\la \X\ra_{n+1}=\la *_{i=1}^{n+1}\la \X\ra\ra.$ On the other hand, from \ref{rkepsilon} (b), 
we have $\varepsilon_n^\wedge(\la \X\ra)=*_{i=0}^{n}\la \X\ra[i]=*_{i=1}^{n+1}\la \X\ra$ since 
$\la \X\ra[i]=\la \X\ra$ for any $i\in\Enteros;$ proving that $\la \X\ra_{n+1}=\la\varepsilon_n^\wedge(\la \X\ra)\ra.$ Similarly, by the dual of \ref{rkepsilon} (b), it follows that $\la \X\ra_{n+1}=\la\varepsilon_n^\vee(\la \X\ra)\ra.$
\

(b) It follows from (a) and \ref{rkepsilon} (c), since $0\in\la \X\ra.$ 
\end{dem}
\vspace{.2cm}
%
%
\vspace{.2cm}

In what follows, we introduce the relative Rouquier's dimension as follows.

\begin{defi} Let $\T$ be a triangulated category, $\X$ a class of objects in $\T$ and $M\in\T.$
The $\X$-dimension of $M$ is $$\dim_\X(M):=\mini\{n\in\N\text{ such that }
M\in\la \X\ra_{n+1}\}.$$
\noindent For a class $\Y$ of objects in $\T,$ we set $\dim_\X(\Y):=\sup\,\{
\dim_\X(Y): \quad\forall\,Y\in\Y\}.$
\end{defi}

\begin{lem}\label{XdimY} Let $\T$ be a triangulated category, and $\X,\,\Y$ be classes of objects in $\T.$ Then, the following statements hold.
\begin{itemize} 
 \item[(a)] For any $n\in\N,$ $\dim_\X(\Y)\leq n$ if and only if $\Y\subseteq\la \X\ra_{n+1}.$
 \item[(b)] If $\X\subseteq\Y$ then $\dim_\Y(M)\leq\dim_\X(M)$ $\;\forall\;M\in\T.$
\end{itemize}
\end{lem}
\begin{dem} (a) Let $\dim_\X(\Y)\leq n.$ Hence, for all $Y\in\Y,$ we have 
$m(Y):=\dim_\X(Y)\leq n.$ Therefore, from \ref{conexion} (b), it follows that 
$Y\in\la \X\ra_{m(Y)+1}\subseteq\la \X\ra_{n+1};$ proving that $\Y\subseteq\la \X\ra_{n+1}.$
Finally, assume that $\Y\subseteq\la \X\ra_{n+1}.$ So, it follows directly that $\dim_\X(\Y)\leq n$
\

(b) Let $\X\subseteq\Y$ and $M\in\T.$ Thus $\la\X\ra\subseteq\la\Y\ra$ and hence 
$*_{i=1}^n\la\X\ra\subseteq *_{i=1}^n\la\Y\ra\subseteq\la *_{i=1}^n\la\Y\ra\ra=\la\Y\ra_n.$ Therefore $\la\X\ra_n\subseteq\la\Y\ra_n$ and so we get that $\dim_\Y(M)\leq\dim_\X(M).$ 
\end{dem}
\vspace{.2cm}

Now, the relative Rouquier's dimension is related to the Rouquier's dimension as follows.

\begin{pro} Let $\T$ be a triangulated category. Then, 
$$\dim\,(\T)=\mini\,\{\dim_X(\T): \quad\forall\; X\in\T\}.$$
\end{pro}
\textbf{Proof}. From \ref{XdimY} (a), we can write down the following equalities
$$\begin{array}{lcl}
\dim\,(\T)&=&\mini\{n\in\N\;\;|\;\text{ there exists }X\in\T\text{ such that }
\la X\ra_{n+1}=\T\}\\
&=&\mini\{n\in\N\;\;|\;\text{ there exists }X\in\T\text{ such that }
\dim_X(\T)\leq n\}\\
&=&\mini\{\dim_X(\T)\quad\forall\; X\in\T\}. \mbox {\hspace{.2cm}  \Qed  }
\end{array}
$$

\vspace{.2cm}

The following are some relationships between relative Rouquier's dimension and the other relative 
dimensions as coresolution, resolution, relative projective and relative injective.

\begin{pro}\label{dimrescores} Let $\T$ be a triangulated category, $\X$ a class of objects in $\T$ and $M\in\T.$
Then, we have that 
$$\dim_\X(M)\leq\maxi\,\{\resdim_{\la \X\ra}(M),\coresdim_{\la \X\ra}(M)\}.$$
\end{pro}
\begin{dem} From \ref{conexion} (a), we have the inclusions $\varepsilon_n^\wedge(\la \X\ra)
\subseteq \la \X\ra_{n+1}$ and $\varepsilon_n^\vee(\la \X\ra)\subseteq \la \X\ra_{n+1}.$ Hence the result follows.
\end{dem}

\begin{pro} Let $(\X,\omega)$ be a pair of classes of objects in $\T$ which are closed under direct summands in $\T.$ If $\X$ is closed under extensions and $\omega$ is an $\X$-injective weak-cogenerator in $\X,$ then
 $$\dim_\X(C)\leq\pd_\omega(C)<\infty,\quad\forall\; C\in\X^\wedge.$$
\end{pro}
\begin{dem} It follows from \ref{dimrescores} and \ref{pd-resd}, since $\X\subseteq \la \X\ra.$
\end{dem}

\begin{pro} Let $(\X,\omega)$ be a pair of classes of objects in $\T$ which are closed under direct summands in $\T.$ If $\X$ is 
cosuspended and $\omega$ is an $\X$-injective weak-cogenerator in $\X,$ then
 $$\dim_\omega(C)\leq\id_\X(C)=\id_\omega(C)<\infty,\quad\forall\; C\in\omega^\vee.$$
\end{pro}
\begin{dem} By \ref{idwX}, we know that $\id_\omega(C)=\id_\X(C)<\infty$ for any $C\in\omega^\sim.$ On the other hand, from \ref{coinyec0} (b),
it follows that $\omega={}_\X\U \cap {}_\X\U^\perp[-1]$ since $\X={}_\X\U.$ In particular, the dual of \ref{pdyresdim} (b) holds (take $\Y:=\omega$)
and hence $\id_\X(C)=\coresdim_\omega(C)$ for any $C\in\omega^\vee.$ Observe that $\omega^{\vee}\subseteq\omega^{\sim}$ since 
$\omega^{\sim}$ is the smallest triangulated subcategory of $\T$ containing $\omega$ (see \ref{tilde-id} (b)). Therefore $\id_\omega(C)=\id_\X(C)<\infty$ for any $C\in\omega^\vee.$ Finally, from
\ref{dimrescores}, we get that $\dim_\omega(C)\leq\coresdim_\omega(C)$ since $\omega\subseteq\la\omega\ra,$ proving the result.
\end{dem}
\vspace{.2cm}

In what follows, $\T$ is assumed to be a $k$-linear, Hom-finite triangulated and Krull-Schmidt 
category over a fixed field $k;$ and $\mathcal{C}$ is a Krull-Schmidt subcategory of $\T$ which 
is closed under direct summands in $\T.$ It is said that $\mathcal{C}$ is $n$-cluster tilting 
(see \cite{KR}) if it is functorially finite and $\mathcal{C}=\cap_{i=1}^{n-1}\,\mathcal{C}[-i]^\perp=\cap_{i=1}^{n-1}\,{}^\perp\mathcal{C}[i].$ 

The following is an example of a triangulated category $\T$ having finite $\mathcal{C}$-dimension.

\begin{pro} Let $\mathcal{C}$ be a $n$-cluster tilting subcategory of $\T.$ Then, 
$$\dim_\mathcal{C}(\T)\leq\resdim_\mathcal{C}(\T)\leq n-1.$$
\end{pro}
\begin{dem} From \cite[Theorem 3.1]{IY}, we know that $\T=*_{i=0}^{n-1}\,\mathcal{C}[i].$ Therefore, by \ref{rkepsilon} (b), we conclude that $\resdim_\mathcal{C}(\T)\leq n-1.$ Finally, by 
\ref{dimrescores}, it follows that $\dim_\mathcal{C}(\T)\leq\resdim_\mathcal{C}(\T)$ since 
$\mathcal{C}\subseteq\la\mathcal{C}\ra;$ proving the result.
\end{dem}
\vspace{.4cm}

\textbf{Acknowledgement} The authors are very grateful to the referee for the comments, corrections and suggestions.

\footnotesize

\vskip3mm \noindent Octavio Mendoza Hern\'andez:\\ Instituto de Matem\'aticas, Universidad Nacional Aut\'onoma de M\'exico\\
Circuito Exterior, Ciudad Universitaria,
C.P. 04510, M\'exico, D.F. MEXICO.\\ {\tt omendoza@matem.unam.mx}

\vskip3mm \noindent Edith Corina S\'aenz Valadez:\\ Departamento de Matem\'aticas,
Facultad de Ciencias, Universidad Nacional Aut\'onoma de M\'exico\\ Circuito Exterior, Ciudad
Universitaria,
C.P. 04510, M\'exico, D.F. MEXICO.\\ {\tt ecsv@lya.fciencias.unam.mx}

\vskip3mm \noindent Valente Santiago Vargas:\\ Instituto de Matem\'aticas, Universidad Nacional Aut\'onoma de M\'exico\\
Circuito Exterior, Ciudad Universitaria,
C.P. 04510, M\'exico, D.F. MEXICO.\\ {\tt valente@matem.unam.mx}

\vskip3mm \noindent Mar\'{i}a Jos\'e Souto Salorio:\\ Facultade de Inform\'atica,
Universidade da Coru\~na\\ 15071 A Coru\~na, ESPA\~NA.\\ {\tt mariaj@udc.es}

\end{document}